\documentclass[11pt]{article}
\usepackage{t1enc}
\usepackage{amsthm}

\newcommand{\N}{\mbox{I\kern-2ptN}}
\newcommand{\R}{\mbox{I\kern-2ptR}}

\newtheorem{teo}{Theorem}
\newtheorem{defi}[teo]{Definition}

\newtheorem{co}[teo]{Corollary}
\newtheorem{lema}[teo]{Lemma}
\theoremstyle{definition}
\newtheorem{nota}[teo]{Remark}

\begin{document}

\title{Minimal strong digraphs}
\author{J. García-López $^{\rm a ,\hspace*{.0cm}}$\thanks{Partially supported by MTM2008-04699-C03-02/MTM.} $\, ,$
 C. Mariju\'an $^{\rm b,\hspace*{-.1cm}}$ \thanks{Partially supported by MTM2007-64704.}$^{\;\; , \hspace*{-.1cm}}$
\thanks{Corresponding author.\newline
{\it E-mail addresses:} jglopez@eui.upm.es (J. García-López), marijuan@mat.uva.es
(C. Mariju\'an).}
}
\date{\small {\it
$^{\rm a}$Dept. Matem\'atica Aplicada, E.U. Inform\'atica, Carretera de Valencia Km 7, 28031-Madrid, Spain\\
$^{\rm b}$Dpto. Matem\'atica Aplicada, E.T.S.I. Inform\'atica, Paseo de Belén 15, 47011-Valladolid, Spain
}}

\maketitle
\begin{abstract}
{We introduce adequate concepts of expansion of a digraph to obtain a sequential construction of minimal
strong digraphs. We characterize the class of minimal strong digraphs whose expansion preserves the property
 of minimality. We prove that every minimal strong digraph of order $n\geq 2$ is the expansion of a
 minimal strong digraph of order $n-1$ and we give sequentially generative procedures for the
constructive characterization of the classes of minimal strong digraphs.
Finally we describe algorithms to compute unlabeled minimal strong digraphs and their isospectral classes}. \\

\noindent {\it 2010-AMS Classification:} 05C20, 05C40, 05C75.
\end{abstract}
\section{Introduction }\label{S1}

In this article, we focus on the study of strongly connected digraphs containing the least possible number of
arcs (minimal strong digraphs), that is, strongly connected digraphs which cease to be so if any one of its arcs
 is suppressed. Minimal strong digraphs can be said to generalize the trees when we consider directed graphs
instead of simply graphs. Nevertheless, the structure of minimal strong digraphs is much richer than that of the trees.

We are previously interested in the following nonnegative inverse eigenvalue problem \cite{To}:
given $k_1, k_2,\dots, k_n$ real numbers, find necessary and sufficient conditions for the existence of a
nonnegative matrix $A$ of order $n$ with characteristic polynomial $x^n+k_1x^{n-1}+ k_2 x^{n-2}+\dots+ k_n$.
The coefficients of the characteristic polynomial are closely related to the cycle structure of the weighted
digraph with adjacency matrix $A$ \cite{CDS}, and the irreducible matricial realizations of the polynomial
are identified with strongly connected digraphs (henceforth strong digraphs) \cite{BR}. The class of strong digraphs
 can easily be reduced to the class of minimal strong digraphs, so we are interested in any theoretical or
constructive characterization of these classes of digraphs.

Many classes of connected graphs and digraphs have constructive characterizations. In particular, for (minimal)
 $2$-connected graphs and (minimal) strong digraphs different procedures have been described to
construct larger (di)graphs from smaller (di)graphs of these classes \cite{Di, Pl, Gr, Ge, He, Bh}.
 The common basic idea of these procedures consists of adding paths between qualified vertices in a systematic way.

Bhogadi \cite{Bh} gives a characterization of Cunningham's decomposition trees for minimal strong digraphs under
X-joint (substitution) composition \cite{Cu}.
He uses his characterization to generate all minimal strong digraphs through 12 vertices and all minimal
$2$-connected graphs through 13 vertices.

As far as we know, these procedures have been defined so that the property of minimality is
 not preserved and the conditions under which minimality is preserved are not characterized. This is not a
difficulty when proving the possibility of obtaining any minimal strong digraph from another smaller one
(Hedetneimi \cite{He} gives a proof by induction), but it is a difficulty when constructing efficient
and sequential procedures and algorithms.

The rest of this paper is organized as follows:

In Section 2, we record basic facts and ideas about the (minimal) strong digraphs.

In Section 3, we introduce two suitable (internal and external) concepts of expansion of a digraph (similar to
the operations ``subdivision'' and ``simple path insertion'' considered by Hedetneimi \cite{He})
for a sequential construction of minimal strong digraphs. We characterize
 the class of minimal strong digraphs whose expansion preserves the property of minimality and we show how every
minimal strong digraph of order $n\geq 2$ is the expansion of a minimal strong digraph of order $n-1$.

In Section 4, we propose a sequentially generative procedure for the constructive characterization of the
 class of minimal strong digraphs.

In Section 5, we implement an algorithm to compute unlabeled minimal strong digraphs following the
 construction of the previous sections. Another algorithm allows the digraphs and the characteristic
 polynomials of the isospectral classes of the minimal strong digraphs to be obtained. \\

\section{Basic general ideas} \label{S3}

In this paper we will use some standard basic concepts and results about graphs as they have been
described in \cite{Ha}.

By a {\bf digraph} $D$ we mean a couple $D=(V,A)$ where $V$ is a finite nonempty set and
$A\subset V\times V-\{(v,v): v\in V\}$.
If $u,v\in V$ we denote $(u,v)$ by $uv$ and we write $D-uv$ and $D+uv$ for the digraphs
$(V,A- \{(u,v)\})$ and $(V,A\cup \{(u,v)\})$, respectively.
For a vertex $v\in V$, the subdigraph $D-v$ consists of all vertices of $D$ except $v$ and all arcs of $D$ except
those incident with $v$.
By a $q$-{\bf cycle} we mean a directed cycle of length $q$ and it is denoted by $C_q$. By a {\bf directed tree}
we mean the digraph obtained from a tree by replacing each edge $\{u,v\}$ with the two arcs $(u,v)$ and $(v,u)$.

A digraph $D$ is {\bf strongly connected} or (simply) {\bf strong} if
every two vertices in $D$ are joined by a path.
It is well known that the digraph $D$ is strongly connected
if and only if its adjacency matrix $M$ is irreducible \cite{BR}.

We record now a number of basic facts about the strong digraphs
that, for simplicity, in the following we will write as SC digraphs. In an SC digraph of order $n\geq 2$
the indegree and outdegree of the vertices are bigger than or equal to $1$.
We will call a vertex with indegree and outdegree equal to $1$ a {\bf linear vertex}.

 If we add an arc to the set of arcs of an SC digraph $D$ then the cyclic structure of $D$ is modified.
This suggests the introduction of the concept of minimal strong digraph.
A SC digraph $D$ is said to be {\bf minimal} if $D-a$ is not strongly connected for any arc $a\in A$.
For simplicity, in the following we will write minimal strong digraph as MSC digraph.

The set of SC digraphs of order $n$ with vertex set $V$ can be partially ordered by the relation of
 inclusion among their sets of arcs. Then, the MSC digraphs are the minimal elements of this partially
 ordered set. Analogously, the set of irreducible $(0,1)$-matrices of order $n$ with zero trace can be
 partially ordered by means of the coordinatewise ordering. The minimal elements of this partially
ordered set are called {\bf nearly reducible} matrices and so the digraph $D$ is an MSD digraph
if and only if its adjacency matrix $M$ is a nearly reducible matrix \cite{BR,HS}. Harfiel \cite{Hf}
gives a remarkably canonical form for nearly reducible matrices.

To reduce the cyclic structure of a SC digraph to the structure of a MSC digraph requires to characterize
the MSC digraphs and to build the set of SC digraphs starting from the set of MSC digraphs.

If $D$ is an MSC digraph and there is a $u-v$ path in $D$, then there cannot be an arc joining the
vertex $u$ to the vertex $v$, i.e. $uv\notin A$. In general, we will say that an arc $uv$ in a digraph $D$
is {\bf transitive} if there is another $u-v$ path distinct from the arc $uv$. We will also say that the
semicycle consisting of a $u-v$ path together with the arc $uv$ is a {\bf pseudocycle}.
So an MSC digraph has no transitive arcs or pseudocycles; moreover, this condition characterizes the minimality
of the strong connection.

\begin{lema} \label{Hedetneimi}
(Geller \cite{Ge}, Hedetneimi \cite{He}) If $D$ is an SC digraph, then $D$ is minimal if and only if $D$ has
no transitive arcs if and only if $D$ has no pseudocycles.
\end{lema}

Consequently, if $D$ is an MSC digraph then so is every strong subdigraph of $D$.

The {\bf contraction} of a cycle in an SC digraph consists of the reduction of the cycle to a unique
vertex, so that $n-1$ of its vertices and its $n$ arcs are eliminated.

\begin{lema} \label{contraction}
(Berge \cite{Be}) The contraction of a cycle in an MSC digraph preserves the minimality,
 i.e. it produces another MSC digraph.
\end{lema}

 The size of an SC digraph of order $n\geq 2$ verifies $n\leq Card(A) \leq n^2-n$ and the
extreme digraphs are the cycle $C_n$ and the complete digraph $K_n$.
The following result was basically obtained by Gupta \cite{Gu}. Brualdi-Hedrick \cite{BHe} gave a different
proof for a more thorough result. We use lemma \ref{contraction} for a shorter proof of the result of Brualdi-Hedrick.

\begin{lema} \label{size}
The size of an MSC digraph $D$ of order $n\geq 2$ verifies $n\leq Card(A)\leq 2(n-1)$.
The size of $D$ is $n$ if and only if $D$ is an $n$-cycle.
The size of $D$ is $2(n-1)$ if and only if $D$ is a directed tree.
\end{lema}

\noindent {\it Proof:}
It is clear that $n\leq Card(A)$ and that the cycle $C_n$ is the unique MSC digraph of order $n$.

Let us see that $Card(A) \leq 2(n-1)$.
We proceede by induction over the order $n$. If $n=2$ the unique MSC digraph is the cycle $C_2$ and the
 inequality is clear for $Card(A)=2$.

\noindent Induction hypothesis: we suppose that every MSC digraph of order $p\leq n$ has at most $2(p-1)$ arcs.

If the MSC digraph is the cycle $C_{n+1}$ the inequality is clear. If $D$ is an MSC digraph of order
$n+1$ distinct from the cycle $C_{n+1}$, as it is an SC digraph, $D$ contains at least a cycle $C_p$ with
$2\leq p\leq n$.
By Lemma \ref{contraction}, the contraction of the cycle $C_p$ produces an MSC digraph $D'$ of order
$n+1-(p-1)=n-p+2\leq n$. By the induction hypothesis, $D'$ has at most $2(n-p+1)$ arcs. Then the number of arcs
of the original digraph $D$ will be at most $2(n-p+1)+p = 2n-p+2\leq 2n$.

Let us see that if $D$ is an MSC digraph of order $n$ and size $2(n-1)$ then it is a directed tree.
Note that the cycles in a directed tree have length two. We suppose, by reductio ad absurdum,
that $D$ has some cycle $C_q$ of length $q>2$. Let $D'$ be the MSC digraph obtained by the
contraction of the cycle $C_q$ in $D$. The order and the size of $D'$ are $n'=n-(q-1)$ and $m'=2(n-1)-q$,
 respectively. Then we have the contradiction $m'\leq 2(n'-1)=2(n-(q-1)-1)=2n-2q<2n-2-q=m'$.
\mbox{}\hfill\rule{.2cm}{.2cm}

Brualdi-Hedrick \cite{BHe} also proved that
there exists an MSC digraph of order $n\geq 2$ and size $m$ if and only if $n\leq m\leq 2(n-1)$
and characterized the MSC digraphs of order $n$ and size $2n-3$.

The next theorem was first proved by Dirac \cite{Di} and independently by Plummer \cite{Pl} in the context of
minimal two connected graphs and by Berge and by Brualdi-Ryser \cite{BR} for minimal strong digraphs.
Our proof is a simplification of that given by Berge \cite{Be}.

\begin{teo} \label{twolv}
 Every MSC digraph of order $n\geq 2$ has at least two linear vertices.
\end{teo}

\noindent {\it Proof:} By induction over the order $n$. If $n=2$ the unique MSC digraph is the cycle $C_2$ whose
vertices are linear.

\noindent Induction hypothesis: we suppose that every MSC digraph of order $p\leq n$ has at least
two linear vertices.

a) If the MSC digraph is the cycle $C_{n+1}$, it has $n+1\geq 3$ linear vertices.

b) If $D$ is an MSC digraph of order $n+1$ that contains no cycle of length bigger than two
then, as it is an SC digraph, it is a directed tree. The extreme vertices (the leaves) of this tree
are the linear vertices of $D$. Because every tree has at least two leaves, then there are at least two
linear vertices in $D$.

c) If $D$ is an MSC digraph of order $n+1$ that contains a cycle $C_p$ of length $p$ with $3\leq p< n+1$,
then there is at least a vertex $v$ in $D$ that is not in the cycle $C_p$. By Lemma \ref{contraction},
the contraction of the cycle $C_p$ produces a new MSC digraph $D'$ of order $n+1-(p-1)=n-p+2$ with $2\leq n-p+2< n$.
By the induction hypothesis, $D'$ has at least two linear vertices that we will call $u$ and $v$. If one of these
vertices, let us suppose the $u$, is the contracted vertex, then in the digraph $D$ there is a unique arc going
into the cycle $C_p$ and a unique arc leaving the cycle $C_p$ and, as $p\geq 3$, in $C_p$ there is at least
one linear vertex $w$. Then $w$ and $v$ are two linear vertices in $D$. If, on the contrary, the linear vertices
$u$ and $v$ of $D'$ are distinct from the contracted vertex, then these vertices are also linear in $D$.
\mbox{}\hfill\rule{.2cm}{.2cm}\\

\section{Sequential expansion of MSC digraphs}

In this section, we look at that any MSC digraph of order $n$ can be generated from an
MSC digraph of order $n-1$. For this purpose, we shall define two different (internal and external) expansion
procedures of a digraph consisting in adding a new vertex so that, either the property of being
 MSC is preserved or the conditions in which the expansion can be carried out while
 preserving the MSC property are described.

The internal expansion ({\it one-step expansion} in \cite{HNC}) of a digraph consists in the sustitution
of an arc $uw$ by new arcs $uv$ and $vw$, $v$ being a new vertex in the digraph. More precisely,

\begin{defi}
The {\bf internal expansion} of the digraph $D=(V,A)$ {\bf by the vertex} $v\notin V$ {\bf over the arc} $uw$
is the digraph $i_{uw}(D)=(V\cup \{v\},A^*)$ with $A^*=A\cup \{uv,vw\}-\{uw\}$.
\end{defi}

The external expansion of a digraph consists in the joining of two vertices $u$ and $w$ (not necessary distinct)
 of the digraph with a new vertex $v$ by means of the arcs $uv$ and $vw$. More precisely,

\begin{defi}
The {\bf external expansion} of the digraph $D=(V,A)$ {\bf by the vertex} $v\notin V$ {\bf from the vertex} $u\in V$
{\bf to the vertex} $w\in V$ is the digraph $e_{uw}(D)=(V\cup \{v\},A^*)$ with $A^*=A\cup \{uv,vw\}$.
\end{defi}

 It is easy to proved that the internal expansion of a digraph preserves the SC and MSC properties
and that the external expansion preserves the SC property but not the MSC property.
The external expansion from the
vertex $u$ to the vertex $w$ can produce transitivity in other arcs,
including when $uw$ is not an arc of an MSC digraph $D$, thus losing the property of minimality.
Next we characterize the class of MSC digraphs whose external expansions preserve the MSC property.

\begin{teo} \label{extexppreservesMSC}
Let $D=(V,A)$ be an MSC digraph and let $u,w$ be vertices such that $uw\notin A$.
The external expansion $e_{uw}(D)$ of $D$ by the vertex $v\notin V$ from the vertex
 $u$ to the vertex $w$ is an MSC digraph if and only if the digraph $D+uw$ has no transitive arcs
distinct from $uw$.
\end{teo}
\noindent {\it Proof:} Clearly $uw$ is a transitive arc of the digraph $D+uw$ because $D$ is an
 SC digraph. If there exists a transitive arc $pq$ distinct from $uw$ in $D+uw$, then there is a longer $p-q$ path
that includes the arc $uw$. This path has the form $p\dots uw\dots q$ where $p$ and $u$ may coincide or
$q$ and $w$ may coincide, but not both simultaneously. Then the path $p\dots uvw\dots q$ makes
 the arc $pq$ transitive in the digraph $e_{uw}(D)$. In fact $pq$ is transitive in $D+uw$ if and only if
$pq$ is transitive in $e_{uw}(D)$ if and only if $e_{uw}(D)$ is not MSC.
\mbox{}\hfill\rule{.2cm}{.2cm}\\

The following result is the base of a possible generative construction of MSC digraphs of order $n\geq 2$
starting from MSC digraphs of order $n-1$. In fact, we shall prove a stronger result;
more exactly, we will prove that every linear vertex of an MSC digraph originate in the (internal
or external) expansion of an MSC digraph. So that if an MSC digraph $D$ has $p\geq 2$ linear vertices,
we can obtain $p$ distinct ``reductions'' with one vertex less than $D$ (some can be isomorph).

\begin{teo}\label{generation}
 Let $D^*=(V,A^*)$ be an MSC digraph of order $n\geq 3$ and $v\in V$ a linear vertex in $D^*$. Then there exists an
 MSC digraph $D=(V-\{v\},A)$ whose (internal or external) expansion by the vertex $v$ is the digraph $D^*$.
\end{teo}

\noindent {\it Proof:} As $v$ is a linear vertex there are two unique vertices $u$ and $w$ such that
$uv\in A^*$ and $vw\in A^*$.

\noindent $a)$ If $u=w$, then $A= A^*-\{uv,vu\}$ and $D=(V-\{v\},A)=D^*-v$ is obviously MSC. By contruction, the
external expansion of the digraph $D$ by the vertex $v$ from the vertex $u$ to the vertex $u$ is the digraph $D^*$.

\noindent $b)$ If $u\neq w$, as there are no transitive arcs in $D^*$, then $uw\notin A^*$.

\noindent $b_1)$ We suppose that no $u-w$ path distinct from the path
$uvw$ exists in $D^*$. In this case we replace the arcs $uv,vw$ in $D^*$ by the new arc $uw$,
 more precisely, we take
$A= A^*\cup \{uw\}-\{uv,vw\}$. The new digraph $D=(V-\{v\},A)$ is by construction SC and,
as there are no $u-w$ paths in $D$, the arc $uw$ is not transitive and then $D$ is also minimal.
By construction, the internal expansion of the digraph $D$ by the vertex $v$ over the arc $uw$
is the digraph $D^*$.

\noindent $b_2)$ If there exists any $u-w$ path distinct from the path $uvw$ in $D^*$, then we make
 $A= A^*-\{uv,vw\}$. The $u-w$ path ensures the strong connection of the new digraph
$D=(V-\{v\},A)=D^*-v$ which is minimal because there are no transitive arcs in $D^*$ and therfore neither in $D$.
 By construction, the external expansion of the digraph $D$ by the vertex $v$ from the vertex $u$ to the
vertex $w$ is the digraph $D^*$.
\mbox{}\hfill\rule{.2cm}{.2cm}\\

\begin{defi}
We will say that the SC digraph $D$ is a {\bf reduction} of the SC digraph $D^*$ if $D^*$ is an internal
or external expansion of $D$.
\end{defi}

 From the above theorems \ref{twolv} and \ref{generation} one can also deduce the following consequences:

\begin{co}
Every MSC digraph of order $\,n\geq 3$ can be reduced to the cycle $C_2$ by a sequence of $\,n-2\,$ reductions.
\end{co}

It is possible to define procedures for the reduction of an MSC digraph to obtain different classes
of MSC digraphs such as a tree $T$ of cycles of distinct lengths, and this tree $T$ can be reduced to one cycle
(whose length is bounded by the biggest of the lengths of the cycles in $T$), or one path of cycles $C_2$ or
one star of cycles $C_2$. All of them can finally be reduced to one cycle $C_2$ and this to a unique vertex.\\

\begin{nota}
Following lemma \ref{contraction}, we can make reductions preserving the MSC property through the contraction
 of cycles. A procedure could be determined by the length of the cycles.
The minimal number of contractions of cycles to reduce an MSC digraph to a vertex is the cyclomatic number
$Card(A)-Card(V)+1$ (Berge, \cite{Be}).
\end{nota}

\section{Construction of MSC and SC digraphs}

In the previous section we saw, on the one hand, that the internal expansion
  of an MSC digraph of order $n$ on any one of its arcs produces an MSC digraph of order $n+1$,
and on the other hand (theorem \ref{extexppreservesMSC}), we saw under which conditions the external expansion
of an MSC digraph of order $n$ over pairs of non adjacent vertices produces an SC digraph of order
 $n+1$ preserving the minimality. We also saw (theorem \ref{generation}) how any MSC digraph of order
$n+1$ can be obtained by (internal or external) expansion of an MSC digraph of order $n$. This three
 results suggests a sequentially generative procedure for the construction of the set of MSC digraphs
 of order $n+1$ starting from the set of MSC digraphs of order $n$. In the figure $1$ we describe the three
 first steps of this process.

\setlength{\unitlength}{0.65cm}
\begin{picture}(9.0,10.0)
\put(0,0){\circle*{.2}}
\put(-0.2,-0.7){\small $u$}

\put(0.5,0){\vector (1,0){1.0}}
\put(0.7,0.3){\small $e_u$}

\put(2,0){\circle*{.2}}
\put(4,0){\circle*{.2}}
\qbezier (2,0) (3.0,0.6) (4,0)
\qbezier (2,0) (3.0,-0.6) (4,-0)
\put(3.5,0.225){\vector (3,-1){0.2}}
\put(2.5,-0.225){\vector (-3,1){0.2}}
\put(1.8,-0.7){\small $u$}
\put(3.8,-0.7){\small $w$}

\put(4.5,0.5){\vector (2,3){1.5}}
\put(4.5,-0.5){\vector (1,-2){1.5}}
\put(4.3,1.7){\small $i_{uw}$}
\put(5.4,-1.8){\small $e_w$}

\put(7,4){\circle*{.2}}
\put(9,4){\circle*{.2}}
\put(8,6){\circle*{.2}}
\put(7,4){\vector (1,2){0.95}}
\put(8,6){\vector (1,-2){0.95}}
\put(9,4){\vector (-1,0){1.9}}
\put(6.8,3.3){\small $u$}
\put(8.8,3.3){\small $w$}
\put(7.9,6.3){\small $v$}

\put(6,-5){\circle*{.2}}
\put(8,-5){\circle*{.2}}
\put(10,-5){\circle*{.2}}
\qbezier (6,-5) (7.0,-4.4) (8,-5)
\qbezier (6,-5) (7.0,-5.6) (8,-5)
\qbezier (8,-5) (9.0,-4.4) (10,-5)
\qbezier (8,-5) (9.0,-5.6) (10,-5)
\put(7.5,-4.775){\vector (3,-1){0.2}}
\put(6.5,-5.225){\vector (-3,1){0.2}}
\put(9.5,-4.775){\vector (3,-1){0.2}}
\put(8.5,-5.225){\vector (-3,1){0.2}}
\put(11.8,-5.7){\small $u$}
\put(13.8,-5.7){\small $w$}
\put(15.8,-5.7){\small $v$}

\put(9.75,4.5){\vector (2,3){1.5}}
\put(9.75,5.7){\small $i_{vw}$}
\put(9.75,4){\vector (1,0){1.5}}
\put(10.2,4.2){\small $e_{w}$}
\put(9.75,3.5){\vector (1,-1){1.5}}
\put(10.6,2.8){\small $e_{wv}$}
\put(10.5,-4.5){\vector (1,1){1}}
\put(10.2,-4.0){\small $i_{uw}$}
\put(10.5,-5){\vector (1,0){1}}
\put(10.8,-4.8){\small $e_{v}$}
\put(10.5,-5.5){\vector (1,-2){1}}
\put(11.1,-6.4){\small $e_{w}$}
\put(12,7.5){\circle*{.2}}
\put(14,7.5){\circle*{.2}}
\put(12,9.5){\circle*{.2}}
\put(14,9.5){\circle*{.2}}
\put(12,7.5){\vector (0,1){1.9}}
\put(12,9.5){\vector (1,0){1.9}}
\put(14,9.5){\vector (0,-1){1.9}}
\put(14,7.5){\vector (-1,0){1.9}}
\put(11.8,6.8){\small $u$}
\put(13.8,6.8){\small $w$}
\put(11.8,9.8){\small $v$}

\put(12,4){\circle*{.2}}
\put(14,4){\circle*{.2}}
\put(13,6){\circle*{.2}}
\put(16,4){\circle*{.2}}
\put(12,4){\vector (1,2){0.95}}
\put(13,6){\vector (1,-2){0.95}}
\put(14,4){\vector (-1,0){1.9}}
\qbezier (14,4) (15.0,4.6) (16,4)
\qbezier (14,4) (15.0,3.4) (16,4)
\put(15.5,4.225){\vector (3,-1){0.2}}
\put(14.5,3.775){\vector (-3,1){0.2}}
\put(11.8,3.3){\small $u$}
\put(13.8,3.3){\small $w$}
\put(12.8,6.3){\small $v$}

\put(12,0.5){\circle*{.2}}
\put(14,0.5){\circle*{.2}}
\put(13,2.5){\circle*{.2}}
\put(15,2.5){\circle*{.2}}
\put(12,0.5){\vector (1,2){0.95}}
\put(13,2.5){\vector (1,-2){0.95}}
\put(14,0.5){\vector (-1,0){1.9}}
\put(14,0.5){\vector (1,2){0.95}}
\put(15,2.5){\vector (-1,0){1.9}}
\put(11.8,-0.2){\small $u$}
\put(13.8,-0.2){\small $w$}
\put(12.8,2.8){\small $v$}

\put(12,-3){\circle*{.2}}
\put(14,-3){\circle*{.2}}
\put(13,-1){\circle*{.2}}
\put(16,-3){\circle*{.2}}
\put(12,-3){\vector (1,2){0.95}}
\put(13,-1){\vector (1,-2){0.95}}
\put(14,-3){\vector (-1,0){1.9}}
\qbezier (14,-3) (15.0,-2.4) (16,-3)
\qbezier (14,-3) (15.0,-3.6) (16,-3)
\put(15.5,-2.775){\vector (3,-1){0.2}}
\put(14.5,-3.225){\vector (-3,1){0.2}}
\put(11.8,-3.7){\small $u$}
\put(13.8,-3.7){\small $w$}
\put(15.8,-3.7){\small $v$}

\put(12,-5){\circle*{.2}}
\put(14,-5){\circle*{.2}}
\put(16,-5){\circle*{.2}}
\put(18,-5){\circle*{.2}}
\qbezier (12,-5) (13.0,-4.4) (14,-5)
\qbezier (12,-5) (13.0,-5.6) (14,-5)
\qbezier (14,-5) (15.0,-4.4) (16,-5)
\qbezier (14,-5) (15.0,-5.6) (16,-5)
\qbezier (16,-5) (17.0,-4.4) (18,-5)
\qbezier (16,-5) (17.0,-5.6) (18,-5)
\put(13.5,-4.775){\vector (3,-1){0.2}}
\put(12.5,-5.225){\vector (-3,1){0.2}}
\put(15.5,-4.775){\vector (3,-1){0.2}}
\put(14.5,-5.225){\vector (-3,1){0.2}}
\put(17.5,-4.775){\vector (3,-1){0.2}}
\put(16.5,-5.225){\vector (-3,1){0.2}}
\put(5.8,-5.7){\small $u$}
\put(7.8,-5.7){\small $w$}
\put(9.8,-5.7){\small $v$}

\put(12,-8.5){\circle*{.2}}
\put(14,-8.5){\circle*{.2}}
\put(16,-8.5){\circle*{.2}}
\put(14,-6.5){\circle*{.2}}
\qbezier (12,-8.5) (13.0,-7.9) (14,-8.5)
\qbezier (12,-8.5) (13.0,-9.1) (14,-8.5)
\qbezier (14,-8.5) (15.0,-7.9) (16,-8.5)
\qbezier (14,-8.5) (15.0,-9.1) (16,-8.5)
\qbezier (14,-8.5) (13.4,-7.5) (14,-6.5)
\qbezier (14,-8.5) (14.6,-7.5) (14,-6.5)
\put(13.5,-8.275){\vector (3,-1){0.2}}
\put(12.5,-8.725){\vector (-3,1){0.2}}
\put(15.5,-8.275){\vector (3,-1){0.2}}
\put(14.5,-8.725){\vector (-3,1){0.2}}
\put(13.75,-7.075){\vector (1,3){0.1}}
\put(14.28,-7.825){\vector (-1,-3){0.1}}
\put(11.8,-9.2){\small $u$}
\put(13.8,-9.2){\small $w$}
\put(15.8,-9.2){\small $v$}

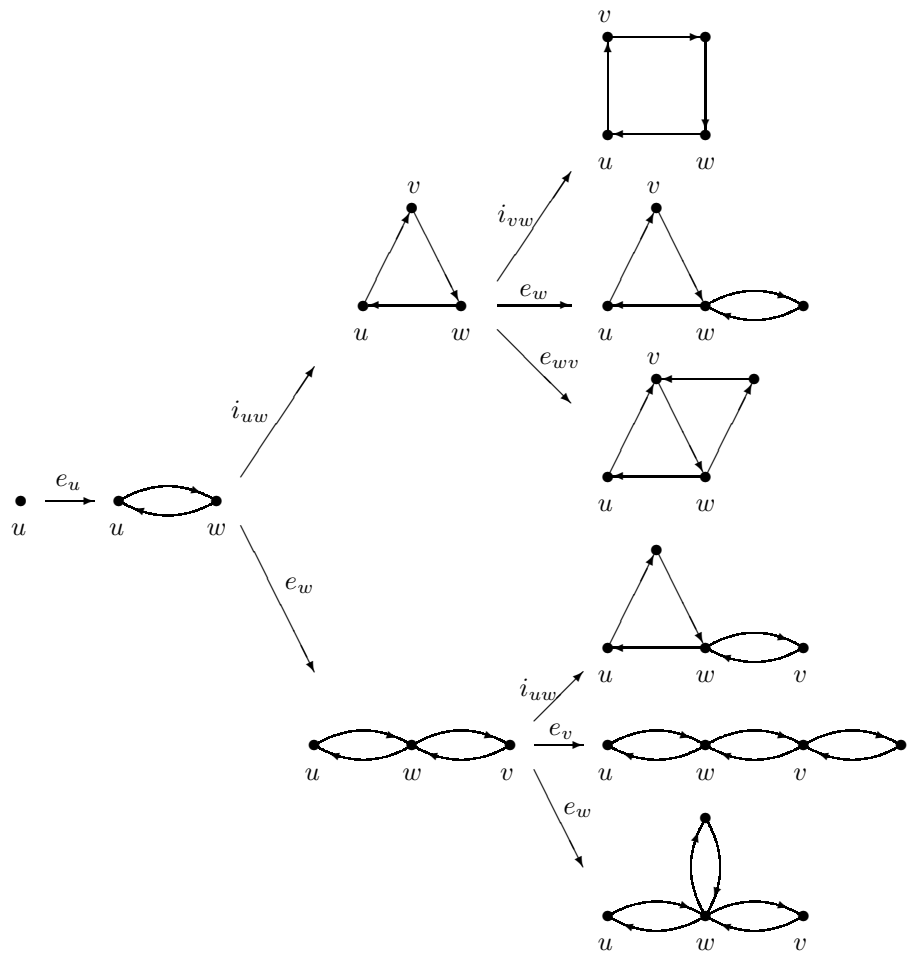
\put(1.8,-10.4){\small Figure $1$. Sequential generative construction of MSC digraphs}

\end{picture}

\vspace{7.3cm}

In general, at the $n$-th iteration, for an MSC digraph $D=(V,A)$ of order $n$ and size $m$ the following is done:

a) an internal expansion over each one of its $m$ arcs;

b) an external expansion over each one of its $n$ vertices;

c) an external expansión from a vertex $u$ to a vertex $w$, such that $uw\notin A$, whenever
(theorem \ref{extexppreservesMSC}) the digraph $D'=(V, A\cup \{uw\})$ has no transitive arcs distinct from $uw$.

 Isomorphism digraphs can be obtained at each step a), b) and c) separately, but also in relation to each other.

To build the set of SC digraphs of order $n$ starting from
 the set of MSC digraphs of order $n$ is sufficient to add any set of transitive arcs.

The above procedures are useful for building and cataloging the sets of MSC digraphs and SC digraphs of order $n$
 but do not give close formulas for the numbers, $UMS(n)$ and $US(n)$, of unlabeled MSC and SC
digraphs of order $n$, respectively.

Labeled strong digraphs were first counted by Liskovec \cite{Li1}, who gives recurrent formulas for the number,
 $S(n)$, of labeled strong digraphs of order $n$ and for the number, $S(n,m)$, of labeled strong
digraphs of order $n$ and size $m$. He also shows the asymptotic behavior $S(n)\approx 2^{n(n-1)}$
 and $US(n)\approx 2^{n(n-1)}/n!$ Liskovec formulas were simplified
by Wright \cite{Wr}, while Robinson \cite{Ro1} gives a natural combinatorial explanation of the simplified equation
of Wright.

Unlabeled strong digraphs were enumerated ``in a somewhat cumbersome manner'' by Liskovec \cite{Li2}
 and Robinson \cite{Ro1} ``outlined'' a method for enumerating them.

The numbers, $MS(n)$ and $UMS(n)$, of labeled and unlabeled MSC digraphs of order $n$ are unknown. \\

\section{Algorithms}

In this section we implement two algorithms. The first one computes unlabeled MSC digraphs,
following the construction described in the previous section. With this algorithm we were able to calculate
all MSC digraphs up to order $14$ in a personal computer. It extends Bhogadi's results to order $13$ and $14$
 and proves the efficiency of our method. Now, we are going to introduce a general description of the algorithm.\\

\noindent {Input:}
\begin{itemize}
\item[] \begin{itemize}
\item[(1)] The order $n$ of the MSC digraphs to be computed.
\item[(2)] The list $L_{n-1}$ of all unlabeled MSC digraphs of order $n-1$.
\end{itemize}
\end{itemize}

\noindent {Output:} A sorted list $L_n$ of all unlabeled MSC digraphs of order $n$.\\

\noindent {Algorithm:}

\begin{itemize}
\item[] \begin{enumerate}
\item[(1)] Set $L=[\ ]$.
\item[(2)] For every $g_{n-1}=(V,A)\in L_{n-1}$:
\begin{itemize}
\item[(a)] For all $uw\in A$:
\begin{itemize}
\item[\bf -] Set $g_n=i_{uw}(g_{n-1})$.
\item[\bf -] Compute $c\_g_n=CanonicalForm(g_n)$
\item[\bf -] If $c\_g_n\not\in L_n$ add the digraph $c\_g_n$ to $L_n$.
\end{itemize}
\item[(b)] For all $u\in V$:
\begin{itemize}
\item[\bf -] Set $g_n=e_{uu}(g_{n-1})$.
\item[\bf -] Compute $c\_g_n=CanonicalForm(g_n)$
\item[\bf -] If $c\_g_n\not\in L_n$ add the digraph $c\_g_n$ to $L_n$.
\end{itemize}
\item[(c)] For all $u\not=w$ such that $uw\not\in A$ and $e_{uw}(g_{n-1})$ is minimal:
\begin{itemize}
\item[\bf -] Set $g_n=e_{uw}(g_{n-1})$.
\item[\bf -] Compute $c\_g_n=CanonicalForm(g_n)$
\item[\bf -] If $c\_g_n\not\in L_n$ add the digraph $c\_g_n$ to $L_n$.
\end{itemize}
\end{itemize}
\end{enumerate}
\end{itemize}

In this algorithm there are three essential procedures. The first one computes a canonical form of a digraph
and it is necessary to detect isomorphic digraphs. Both procedures can be solved by using the software package
{\it nauty} \cite{Kay}. However, for MSC digraphs, we can consider another efficient method. Compute a vertex
 set partition $V_1,\,V_2,\,\dots,\,V_k$ in such a way that, given two arbitrary subsets $V_i$ and $V_j$,
all vertices of $V_i$ have the same number of arcs with the end vertex in $V_j$. Finally, obtain a canonical
 form from this partition. If the canonical form computing has complexity $O(f(n))$ then the overall complexity
of this procedure is $O(n^2|L_{n-1}|f(n))$.

Let $D=(V,A)$ be an MSC digraph and let $u$, $w$ be vertices such that $uw\not\in A$. The second procedure
determines wheter the external expansion $e_{uw}(D)$ is minimal, by using the characterization of theorem
\ref{extexppreservesMSC}.
 For every arc $xz\in D+uw$, $xz\neq uw$, we have to compute whether $xz$ is transitive. Each case can be solved
in $O(n)$ time, checking if there is a path from $x$ to $z$ in the digraph $(D+uw)-xz$. Thus, this procedure has
 complexity $O(n^2)$ and, considering all cases, the overall complexity is $O(n^3|L_{n-1}|)$.

The last procedure updates the sorted list of digraphs $L_n$. It is a very known problem that can be solved
in logarithmic time. However, the size of the list increases very quickly. Therefore, it is necessary to store
the list on a hard disk. Then the overall complexity of this procedure is
$O(n^2|L_{n-1}|\log(n^2|L_{n-1}|))$ because there are $O(n^2|L_{n-1}|)$ updates.

We summarize the results of the computation in table $1$. For every $n$ from $1$ to $14$, it includes the total
 number, $UMS(n)$, of unlabeled MSC digraphs of order $n$. We also classify the MSC digraphs of a given order by the number
 of arcs. When the number of arcs is equal to $2n-2$ the digraphs become directed trees, changing $n$,
the following sequence of unlabeled trees is obtained:
$1,\,1,\,2,\,3,\,6,\,11,\,23,\,47,\,106,\,235,\,551,\,1301,\,3159\dots$.\\

{\tiny
$$
\begin{array}{|c|r|r|r|r|r|r|r|r|r|r|r|r|r|}
\noalign{\hrule}
m\backslash n & 2 & 3 & 4 & 5 & 6 & 7 & 8 & 9 & 10 & 11 & 12 & 13 & 14 \\
\noalign{\hrule}
2 & 1 & & & & & & & & & & & & \\
\noalign{\hrule}
3 & & 1 & & & & & & & & & & & \\
\noalign{\hrule}
4 & & 1 & 1 & & & & & & & & & & \\
\noalign{\hrule}
5 & & & 2 & 1 & & & & & & & & & \\
\noalign{\hrule}
6 & & & 2 & 4 & 1 & & & & & & & & \\
\noalign{\hrule}
7 & & & & 7 & 6 & 1 & & & & & & & \\
\noalign{\hrule}
8 & & & & 3 & 27 & 9 & 1 & & & & & & \\
\noalign{\hrule}
9 & & & & & 23 & 70 & 12 & 1 & & & & & \\
\noalign{\hrule}
10 & & & & & 6 & 131 & 169 & 16 & 1 & & & & \\
\noalign{\hrule}
11 & & & & & & 66 & 559 & 344 & 20 & 1 & & & \\
\noalign{\hrule}
12 & & & & & & 11 & 571 & 1970 & 662 & 25 & 1 & & \\
\noalign{\hrule}
13 & & & & & & & 191 & 3479 & 5874 & 1159 & 30 & 1 & \\
\noalign{\hrule}
14 & & & & & & & 23 & 2229 & 17109 & 15526 & 1947 & 36 & 1 \\
\noalign{\hrule}
15 & & & & & & & & 541 & 18509 & 69845 & 37072 & 3086 & 42 \\
\noalign{\hrule}
16 & & & & & & & & 47 & 8226 & 120582 & 246971 & 81561 & 4743 \\
\noalign{\hrule}
17 & & & & & & & & & 1514 & 87963 & 646339 & 773413 & 167500 \\
\noalign{\hrule}
18 & & & & & & & & & 106 & 28879 & 732150 & 2954946 & 2191491 \\
\noalign{\hrule}
19 & & & & & & & & & & 4217 & 385484 & 4974754 & 11819034 \\
\noalign{\hrule}
20 & & & & & & & & & & 235 & 98146 & 3973379 & 28600421 \\
\noalign{\hrule}
21 & & & & & & & & & & & 11724 & 1587924 & 33313635 \\
\noalign{\hrule}
22 & & & & & & & & & & & 551 & 324638 & 19785730 \\
\noalign{\hrule}
23 & & & & & & & & & & & & 32527 & 6234794 \\
\noalign{\hrule}
24 & & & & & & & & & & & & 1301 & 1052874 \\
\noalign{\hrule}
25 & & & & & & & & & & & & & 90285 \\
\noalign{\hrule}
26 & & & & & & & & & & & & & 3159 \\
\noalign{\hrule}
UMS(n) & 1 & 2 & 5 & 15 & 63 & 288 & 1526 & 8627 & 52021 & 328432 & 2160415 & 14707566 & 103263709 \\
\noalign{\hrule}
\end{array}
$$
}

\centerline{Table 1. Number of unlabeled MSC digraphs of order $n$ and $m$ arcs.}

\vskip0.3cm

The other implemented algorithm computes the isospectral classes of the MSC digraphs. It determines the
 digraphs and the characteristic polynomial of each class. If Gauss's algorithm is used in order to compute
 characteristic polynomials, the overall complexity is $O(n^3|L_n|)$. Table $2$ includes the obtained results.
 Observe that, for $n\geq 8$, there are isospectral classes realized by MSC digraphs with a different number
 of arcs. In order to explain this fact, we have included three summary rows. The first one is the sum of the
numbers of the isospectral classes in the number of possible arcs, the second one includes the total number of
 isospectral classes of a given order and the last one is the difference between them. \\

{\tiny
$$
\begin{array}{|c|r|r|r|r|r|r|r|r|r|r|r|r|r|}
\noalign{\hrule}
m\backslash n & 2 & 3 & 4 & 5 & 6 & 7 & 8 & 9 & 10 & 11 & 12 & 13 & 14 \\
\noalign{\hrule}
2 & 1 & & & & & & & & & & & & \\
\noalign{\hrule}
3 & & 1 & & & & & & & & & & & \\
\noalign{\hrule}
4 & & 1 & 1 & & & & & & & & & & \\
\noalign{\hrule}
5 & & & 2 & 1 & & & & & & & & & \\
\noalign{\hrule}
6 & & & 2 & 4 & 1 & & & & & & & & \\
\noalign{\hrule}
7 & & & & 6 & 6 & 1 & & & & & & & \\
\noalign{\hrule}
8 & & & & 3 & 18 & 9 & 1 & & & & & & \\
\noalign{\hrule}
9 & & & & & 16 & 35 & 12 & 1 & & & & & \\
\noalign{\hrule}
10 & & & & & 6 & 62 & 65 & 16 & 1 & & & & \\
\noalign{\hrule}
11 & & & & & & 43 & 172 & 103 & 20 & 1 & & & \\
\noalign{\hrule}
12 & & & & & & 11 & 227 & 395 & 160 & 25 & 1 & & \\
\noalign{\hrule}
13 & & & & & & & 115 & 801 & 791 & 227 & 30 & 1 & \\
\noalign{\hrule}
14 & & & & & & & 22 & 769 & 2290 & 1423 & 319 & 36 & 1 \\
\noalign{\hrule}
15 & & & & & & & & 319 & 3530 & 5567 & 2411 & 424 & 42 \\
\noalign{\hrule}
16 & & & & & & & & 42 & 2645 & 12437 & 11942 & 3807 & 559 \\
\noalign{\hrule}
17 & & & & & & & & & 848 & 14978 & 36638 & 23583 & 5805 \\
\noalign{\hrule}
18 & & & & & & & & & 102 & 8812 & 64337 & 93732 & 43070 \\
\noalign{\hrule}
19 & & & & & & & & & & 2349 & 61376 & 228358 & 217303 \\
\noalign{\hrule}
20 & & & & & & & & & & 204 & 29317 & 318654 & 695323 \\
\noalign{\hrule}
21 & & & & & & & & & & & 6401 & 244989 & 1351485 \\
\noalign{\hrule}
22 & & & & & & & & & & & 488 & 95369 & 1517405 \\
\noalign{\hrule}
23 & & & & & & & & & & & & 17660 & 949476 \\
\noalign{\hrule}
24 & & & & & & & & & & & & 1078 & 307783 \\
\noalign{\hrule}
25 & & & & & & & & & & & & & 48567 \\
\noalign{\hrule}
26 & & & & & & & & & & & & & 2723 \\
\noalign{\hrule}
sum & 1 & 2 & 5 & 14 & 47 & 161 & 614 & 2446 & 10387 & 46023 & 213260 & 1027691 & 5139542 \\
\noalign{\hrule}
total & 1 & 2 & 5 & 14 & 47 & 161 & 604 & 2360 & 9796 & 42510 & 193891 & 922109 & 4560898 \\
\noalign{\hrule}
\Delta & 0 & 0 & 0 & 0 & 0 & 0 & 10 & 86 & 591 & 3513 & 19369 & 105582 & 578644 \\
\noalign{\hrule}
\end{array}
$$
}

\centerline{Table 2. Isospectral classes of MSC digraphs of order $n$ and $m$ arcs.}

\vskip0.3cm

Finally, we remark that, from this table, we can extract the following sequences of isospectral classes:
\begin{enumerate}
\item[1.] For MSC digraphs: $1,\,2,\,5,\,14,\,47,\,161,\,604,\,2360,\,9796,\,42510,\,193891,$ $922109,\,4560898\dots$.
\item[2.] For trees: $1,\,1,\,2,\,3,\,6,\,11,\,22,\,42,\,102,\,204,\,488,\,1078,\,2723\dots$.
\end{enumerate}

\begin{nota}
With respect to our initial motivation of the nonnegative inverse eigenvalue problem, in the context of
(minimal) strong digraphs we must ask ourselves:

\noindent a) which monic polynomials of degree $n$ with integral coefficients are the characteristic
polynomials of SC digraphs of order $n$, that is to say, of irreducible $(0,1)$-matrices of order $n$
with zero trace, and, analogously

\noindent b) which monic polynomials of degree $n$ with integral coefficients are the characteristic
polynomials of MSC digraphs of order $n$, that is to say, of nearly reducible $(0,1)$-matrices of order $n$.

The problem a) is open and we have nothing relevant to say about it.

The problem b) has been indirectly solved in this paper in the sense that the above algorithms allow the
 class of characteristic polynomials of the nearly reducible matrices of order $n$ and the sets of
 MSC digraphs with equal characteristic polynomial to be catalogued.

The Figure $2$ shows the first pair of non-isomorphic MSC digraphs having the same characteristic polynomial,
in this case $x^5-x^3-2x^2$.

\begin{center}
\setlength{\unitlength}{1cm}
\begin{picture}(8.0,2.5)

\put(0,1){\circle*{.2}}
\put(1,0){\circle*{.2}}
\put(1,2){\circle*{.2}}
\put(2,1){\circle*{.2}}
\put(3,2){\circle*{.2}}

\put(1,0){\vector(-1,1){0.93}}
\put(0,1){\vector(1,1){0.93}}
\put(1,2){\vector(0,-1){1.88}}
\put(1,0){\vector(1,1){0.93}}
\put(2,1){\line(-1,1){0.95}}
\put(1.5,1.5){\vector(-1,1){0.2}}

\qbezier (1,2) (2.0,2.6) (3,2)
\qbezier (1,2) (2.0,1.4) (3,2)
\put(2.5,2.225){\vector (3,-1){0.2}}
\put(1.5,1.775){\vector (-3,1){0.2}}

\put(5,1){\circle*{.2}}
\put(6,0){\circle*{.2}}
\put(6,2){\circle*{.2}}
\put(7,1){\circle*{.2}}
\put(8,0){\circle*{.2}}

\put(6,0){\vector(-1,1){0.93}}
\put(5,1){\vector(1,1){0.93}}
\put(6,2){\vector(0,-1){1.88}}
\put(6,0){\vector(1,1){0.93}}
\put(7,1){\line(-1,1){0.95}}
\put(6.5,1.5){\vector(-1,1){0.2}}

\qbezier (6,0) (7.0,0.6) (8,0)
\qbezier (6,0) (7.0,-0.6) (8,0)
\put(7.5,0.225){\vector (3,-1){0.2}}
\put(6.5,-0.225){\vector (-3,1){0.2}}

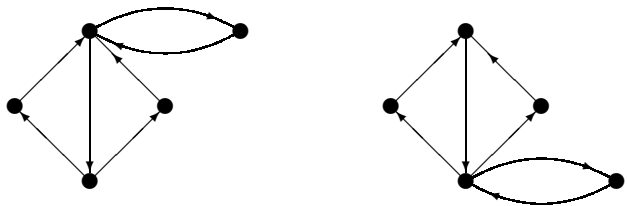
\put(0.0,-0.9){\small Figure 2. Non-isomorphic isospectral MSC digraphs}
\end{picture}
\end{center}

\vskip0.8cm

It is well known that there exist classes of isospectral trees which are as large as desired \cite{CDS}.
So, classes of MSC digraphs (in particular directed trees) can be also be built which can be any size with
the same characteristic polynomial.

It is also well known that the isospectrality relationship does not preserve the connectivity of graphs
 \cite{CDS}. Only the first of the SC digraphs of the figure $3$ is minimal but both have equal characteristic
polynomial $x^5-3x^2$, so the isospectrality relationship does not preserve the minimality of the strong connection
either.
\end{nota}

\begin{center}
\setlength{\unitlength}{1cm}
\begin{picture}(8.0,2.5)

\put(0,1){\circle*{.2}}
\put(1,0){\circle*{.2}}
\put(1,2){\circle*{.2}}
\put(2,1){\circle*{.2}}
\put(3,1){\circle*{.2}}

\put(1,0){\vector(-1,1){0.93}}
\put(0,1){\vector(1,1){0.93}}
\put(1,2){\vector(0,-1){1.88}}
\put(1,0){\vector(1,1){0.93}}
\put(2,1){\line(-1,1){0.95}}
\put(1.5,1.5){\vector(-1,1){0.2}}
\put(1,0){\vector(2,1){1.88}}
\put(3,1){\line(-2,1){1.9}}
\put(1.5,1.75){\vector(-2,1){0.2}}

\put(5,1){\circle*{.2}}
\put(6,0){\circle*{.2}}
\put(6,2){\circle*{.2}}
\put(7,1){\circle*{.2}}
\put(8,2){\circle*{.2}}

\put(6,0){\vector(-1,1){0.93}}
\put(5,1){\vector(1,1){0.93}}
\put(6,2){\vector(0,-1){1.88}}
\put(6,0){\vector(1,1){0.93}}
\put(7,1){\vector(-1,1){0.93}}
\put(6,2){\vector(1,0){1.88}}
\put(8,2){\vector(-1,-1){0.93}}

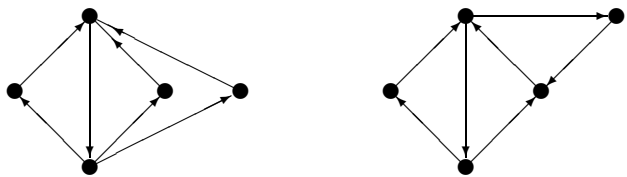
\put(0.6,-0.8){\small Figure 3. MSC and SC isospectral digraphs}
\end{picture}
\end{center}

\end{document}